\documentclass[12pt]{article}

\newcommand{\qed}{{\hfill\rule{4pt}{7pt}}}

\renewcommand{\leq}{\leqslant}
\renewcommand{\geq}{\geqslant}

\usepackage{amsthm,amsmath,amssymb,latexsym,url}

\newtheorem{thm}{Theorem}[section]

\newtheorem{cor}[thm]{Corollary}

\newtheorem{lem}[thm]{Lemma}
\newtheorem{conj}[thm]{Conjecture}

\theoremstyle{remark}
\newtheorem*{rmk}{Remark}

\def\N{\mathbb{N}}

\def\pf{\noindent {\it Proof.} }

\numberwithin{equation}{section}
\renewcommand{\qed}{{\hfill\rule{4pt}{7pt}}\medskip}

\begin{document}

\begin{center}
{\Large\bf Factors of binomial sums from  the Catalan triangle}
\end{center}
\vskip 2mm \centerline{Victor J. W. Guo$^1$  and Jiang Zeng$^{2}$}
\begin{center}
{\footnotesize $^1$Department of Mathematics, East China Normal University,\\ Shanghai 200062,
 People's Republic of China\\
{\tt jwguo@math.ecnu.edu.cn,\quad http://math.ecnu.edu.cn/\textasciitilde{jwguo}}\\[10pt]
$^2$Universit\'e de Lyon; Universit\'e Lyon 1; Institut Camille
Jordan, UMR 5208 du CNRS;\\ 43, boulevard du 11 novembre 1918,
F-69622 Villeurbanne Cedex, France\\
{\tt zeng@math.univ-lyon1.fr,\quad
http://math.univ-lyon1.fr/\textasciitilde{zeng}} }
\end{center}


\vskip 0.7cm {\small \noindent{\bf Abstract.}
By using the Newton interpolation formula, we generalize the recent identities on the Catalan triangle
obtained by Miana and Romero as well as those of Chen and Chu.
We further study  divisibility properties of sums of products of
binomial coefficients and an odd power of a natural number. For example, we prove that
for all positive integers
$n_1, \ldots, n_m$, $n_{m+1}=n_1$, and any nonnegative integer $r$, the expression
$$n_1^{-1}{n_1+n_{m}\choose n_1}^{-1}
\sum_{k=1}^{n_1}k^{2r+1}\prod_{i=1}^{m} {n_i+n_{i+1}\choose n_i+k}$$
is either an integer or a half-integer.
Moreover, several related conjectures are proposed.

}

\vskip 0.2cm
\noindent{\it Keywords:} Catalan triangle; divisibility; Chu-Vandermonde formula; Pfaff-Saalsch\"utz identity;
Lucas' theorem

\vskip 0.2cm
\noindent{\it AMS Subject Classifications} (2000): 05A10, 05A30, 11B65.

\section{Introduction}

Shapiro~\cite{Shapiro}  introduced the  \emph{Catalan triangle} $(B_{n,k})_{n\geq k\geq 0}$, where
$B_{n,k}:=\frac{k}{n}{2n\choose n-k}$ ($1\leq k\leq n$) and
proved that
$\sum_{k=1}^nB_{n,k}=\frac{1}{2}{2n\choose n}$, i.e.,
\begin{align}
\sum_{k=1}^n k{2n\choose n-k}&=\frac{n}{2}{2n\choose n}. \label{eq:shap}
\end{align}
Recently, Guti\'errez et al. \cite{GHMR}, Miana and Romero
\cite{MR}, and Chen and Chu \cite{CC} studied the binomial sums
$\sum_{k=1}^n k^m {2n\choose n-k}^2$. On the other hand, Miana and
Romero \cite{MR} have  proved  the following identity:
\begin{align}\label{eq:mr}
\sum_{p=1}^iB_{n,k}B_{n,n+k-i}(n+2k-i)^3={2n\choose n}{2n-2\choose i-1}(n^2+4n-2ni+i^2),
\end{align}
and asked: Are there polynomials $P_m(n,i)$ and $Q_m(n,i)$ of integral coefficients such that
\begin{align}\label{eq:question}
\sum_{k=1}^iB_{n,k}B_{n,n+k-i}(n+2p-i)^{2m+1}={2n\choose n}{2n-2\choose i-1} \frac{P_m(n,i)}{Q_m(n,i)}
\end{align}
for $m\in \N$ and $1\leq i\leq n$ ?

Our first aim is to give a positive answer to their question. To this end,
consider the sum
\begin{align}\label{eq:deftheta}
\Theta_{2m+1}(n,r)&:=\sum_{\ell=1}^{n-r} \ell^m(\ell+r)^m(2\ell+r) {2n\choose n-\ell}{2n\choose n+\ell+r},
\end{align}
and define the $\alpha$-coefficient  by
\begin{align}
\alpha_k(m,n,r):= \sum_{i=0}^k{2n+r\choose i}{2k-2n-r\choose
k-i}\frac{\big((n-i)^2+r(n-i)\big)^m(2n+r-2i)}{(2n+r-2k)_{2k+1}},
\label{eq:alpha}
\end{align}
where  $(a)_n=a(a+1)\cdots(a+n-1)$ is the Pochhammer symbol.
\begin{thm}\label{thm:1.1}
For $m,n,r\geq 0$, there holds
$$
\Theta_{2m+1}(n,r)
=n{2n\choose n}{2n\choose n-r-1}
\sum_{k=0}^m {n-1\choose k}{n+r+1\choose k+1}\frac{k!(k+1)!}{2n-k}\alpha_k(m,n,r).
$$
\end{thm}
\medskip

\begin{rmk} When  $r=0$ Theorem 1.1 reduces to Theorem 3 in \cite{CC}.
Writing $(2k+r)^{2m}=(4k(k+r)+r^2)^m=\sum_{i=0}^m{m\choose i}4^ik^i(k+r)^i r^{2m-2i}$,
we see from \eqref{eq:deftheta} that
\begin{align}
\sum_{k=1}^{n-r}B_{n,k}B_{n,k+r} (2k+r)^{2m+1}=\sum_{i=0}^m{m\choose i}4^i \Theta_{2i+1}(n,r) r^{2m-2i}, \label{eq:odd}
\end{align}
which  gives a {\it positive} answer to the question \eqref{eq:question} by Theorem~\ref{thm:1.1}.
By the way, according to \cite[(13g)]{CC} we have
$$
\frac{\sum_{p=1}^np^{9}(B_{n,p})^2}{(n+1)C_n  C_{n-2}}=
\frac{n}{2n-5}(30n^5-150n^4+252n^3-185n^2+65n-9).
$$
So the answer to question (2) in \cite{MR} is {\it negative}.
\end{rmk}

Secondly, motivated by the divisibility results in \cite{Calkin,CD,CC,GJZ,GHMR,MR,Shapiro,Zu}, we shall prove the following theorems.

\begin{thm}\label{thm:half}
For all nonnegative integers $n$ and $r\geq 1$,
\begin{align*}
n^{-2}{2n\choose n}^{-1}\sum_{k=1}^n{2n\choose n-k}k^{2r+1}
\end{align*}
is always a half-integer.
\end{thm}

\begin{thm}\label{thm:nnnhalf}
For all positive integers $n_1,\ldots,n_{m}$ and any nonnegative integer $r$,
\begin{align*}
S_{2r+1}(n_1,\ldots,n_{m}):=n_1^{-1}{n_1+n_{m}\choose n_1}^{-1}
\sum_{k=1}^{n_1}k^{2r+1}\prod_{i=1}^{m} {n_i+n_{i+1}\choose n_i+k}\quad (n_{m+1}=n_1)
\end{align*}
is either an integer or a half-integer.
\end{thm}
For example, we have
$$S_7(3,3,2,3,3)=\frac{10233}{2},\qquad
S_{11}(3,3,2,2,2)=2448.
$$

\begin{thm}\label{thm:ppower}
Let $n$ be a prime power. Let $r\geq 0$ and $s\geq 1$ such that $r\not\equiv s\pmod 2$. Then
\begin{align*}
{2n\choose n}^{-1}\sum_{k=1}^{n}k^{r}B_{n,k}^{s}
\end{align*}
is either an integer or a half-integer.
\end{thm}

\begin{rmk}
It seems that
Theorem {\rm\ref{thm:ppower}} holds for any positive integer $n$. For example, one can easily check that
\begin{align*}
\sum_{k=1}^{6}k^{r}B_{6,k}^{s}=132^s+2^r 165^s+3^r 110^s+4^r 44^s+5^r 10^s+6^r
\end{align*}
is divisible by $\frac{1}{2}{12\choose 6}=2\times3\times7\times11$.
\end{rmk}
The rest of this paper is organized as follows. We prove Theorems
1.1, 1.2 , 1.3  and 1.4 in Sections 2, 3, 4 and 5, respectively.
Some consequences of Theorem 1.2 are given in Section 6. Some open
problems and conjectures are proposed in Section 7.

\section{Proof of Theorem~\ref{thm:1.1}}
Let $x_0, x_1, \ldots, x_m$ be different  values of a variable $x$.
Then
the Newton interpolation formula
 (see \cite[Chapter 1]{MT}) implies that
\begin{align}
x^m=\sum_{k=0}^m \left(\sum_{i=0}^k\frac{x_i^m}{\prod_{j=0;\,j\neq i}^k(x_i-x_j)}\right)
\prod_{i=0}^{k-1}(x-x_i). \label{eq:newton}
\end{align}
Letting $x_i=(n-i)^2+r(n-i)$ in \eqref{eq:newton}, we have
\begin{align}
x^m=\sum_{k=0}^m \alpha_k(m,n,r)
\prod_{i=0}^{k-1}\big((n-i)^2+r(n-i)-x\big),  \label{eq:nnrn}
\end{align}
where $\alpha_k(m,n,r)$ is given by \eqref{eq:alpha}.

Replacing $x$ by $x(x+r)$ in \eqref{eq:nnrn} and noticing that
$(n-i)^2+r(n-i)+x(x+r)=(n-x-i)(n+x+r-i)$,
we get
\begin{align}
x^m(x+r)^m=\sum_{k=0}^m {n-x\choose k}{n+x+r\choose k}k!^2\alpha_k(m,n,r).  \label{eq:ellr}
\end{align}
Using \eqref{eq:ellr} to rewrite  $\ell^m(\ell+r)^m$ and
applying the binomial relation
$$
{n-\ell\choose k}{n+\ell+r\choose k}{2n\choose n-\ell}{2n\choose n+\ell+r}
={2n\choose k}^2{2n-k\choose n+\ell}{2n-k\choose n-\ell-r},
$$
one has
\begin{align}
\Theta_{2m+1}(n,r)
&=\sum_{\ell=1}^{n-r}  (2\ell+r) {2n\choose n-\ell}{2n\choose n+\ell+r}
\sum_{k=0}^m {n-\ell\choose k}{n+\ell+r\choose k}\alpha_k(m,n,r)\nonumber\\
&=\sum_{k=0}^m {2n\choose k}^2\alpha_k(m,n,r)
\sum_{\ell=1}^{n-r} (2\ell+r) {2n-k\choose n+\ell}{2n-k\choose n-\ell-r}.  \label{eq:theta}
\end{align}
Noticing that
\begin{align*}
&\hskip -3mm (2\ell+r) {x\choose n+\ell}{x\choose n-\ell-r}\\
&= (n+\ell){x\choose n+\ell}{x-1\choose n-\ell-r}
-(n+\ell+1){x\choose n+\ell+1}{x-1\choose n-\ell-r-1},
\end{align*}
we have
\begin{align*}
\sum_{\ell=1}^{n-r} (2\ell+r) {2n-k\choose n+\ell}{2n-k\choose n-\ell-r}
=(n+1){2n-k\choose n+1}{2n-k-1\choose n-r-1}.
\end{align*}
The theorem then follows by observing that
$$
(n+1){2n\choose k}^2{2n-k\choose n+1}{2n-k-1\choose n-r-1}
=n{2n\choose n}{2n\choose n-r-1}{n-1\choose k}{n+r+1\choose k+1}
\frac{k+1}{2n-k}.
$$

For the reader's convenience, we derive the first values of $\Theta_{2m+1}(n,r)$ from Theorem~\ref{thm:1.1}.
\begin{align*}
\Theta_1(n,r)&={2n\choose n}{2n-2\choose n-r-1}{n^3},\\[5pt]
\Theta_3(n,r)&={2n\choose n}{2n-2\choose n-r-1}\frac{n^3(3n^2-5n-r^2+1)}{2n-3},\\[5pt]
\Theta_5(n,r)&={2n\choose n}{2n-2\choose n-r-1}\frac{n^3(6n^3-12n^2-4nr^2+6n+r^2-1)}{2n-3},\\[5pt]
\Theta_7(n,r)&=\frac{{2n\choose n}{2n-2\choose n-r-1}n^3}{(2n-3)(2n-5)}
\left\{
\begin{array}{l}
30n^5-150n^4+252n^3-30n^3r^2+91n^2r^2-185n^2\\
{}-53nr^2+4nr^4+65n-r^4+10r^2-9
\end{array}\right\}.
\end{align*}
Therefore, by \eqref{eq:odd} we have
\begin{align}
&\sum_{k=1}^{n-r} B_{n,k}B_{n,k+r}(2k+r)^3={2n\choose n}{2n-2\choose n-r-1}(4n+r^2), \label{eq:newkr} \\
&\sum_{k=1}^{n-r} B_{n,k}B_{n,k+r}(2k+r)^5={2n\choose n}{2n-2\choose n-r-1}\left(\frac{16n(3n^2-5n-r^2+1)}{2n-3}+8nr^2+r^4\right), \nonumber \\
&\sum_{k=1}^{n-r} B_{n,k}B_{n,k+r}(2k+r)^7={2n\choose n}{2n-2\choose n-r-1}  \nonumber\\
&\quad{}\times\left(\frac{64n(6n^3-12n^2-4nr^2+6n+r^2-1)+48n(3n^2-5n-r^2+1)r^2}{2n-3}+12nr^4+r^6\right).  \nonumber
\end{align}
Eq. \eqref{eq:newkr} is an equivalent form of  \cite[Theorem 2.3]{MR}, i.e., \eqref{eq:mr}.

\section{Proof of Theorem \ref{thm:half}}

\pf
As in \eqref{eq:nnrn} we can prove
the following identity (see \cite[(9)]{CC}) by the Newton interpolation formula:
\begin{align}
\ell^{2r}=\sum_{k=0}^r {n-\ell\choose k}{n+\ell\choose k}
\frac{2k!^2}{(2n-2k)_{2k+1}}\sum_{i=0}^k{2n\choose i}{2k-2n\choose k-i}(n-i)^{2r+1}. \label{eq:ellnor}
\end{align}
It follows from \eqref{eq:ellnor} and \eqref{eq:shap} that
\begin{align}
\sum_{\ell=1}^n{2n\choose n-\ell}\ell^{2r+1}
&=\sum_{\ell=1}^n\sum_{k=0}^r {2n-2k\choose n-k-\ell}
\frac{\ell}{n-k}\sum_{i=0}^k{2n\choose i}{2k-2n\choose k-i}(n-i)^{2r+1} \nonumber\\
&=\frac{1}{2} \sum_{k=0}^rf_{n,k}(r),  \label{eq:sumnkl}
\end{align}
where
\begin{align*}
f_{n,k}(r):={2n-2k\choose n-k}
\sum_{i=0}^k{2n\choose i}{2k-2n\choose k-i}(n-i)^{2r+1}.
\end{align*}
We now show by   induction on $r$ that  $f_{n,k}(r)$
is divisible by $n^{\min\{2,r+1\}}{2n\choose n}$. For $r=0$,  writing $f_{n,k}(0)$
as
\begin{align*}
2n{2n-2k\choose n-k}
\sum_{i=0}^k{2n-1\choose i}{2k-2n\choose k-i}
-n{2n-2k\choose n-k}
\sum_{i=0}^k{2n\choose i}{2k-2n\choose k-i},
\end{align*}
we see, by the Chu-Vandermonde formula, that
\begin{align*}
f_{n,k}(0)=n{2n-2k\choose n-k}\left(2{2k-1\choose k}-{2k\choose k}\right)
=\begin{cases}
0, &\text{if $k>0$}, \\[5pt]
\displaystyle n{2n\choose n}, &\text{if $k=0$}.
\end{cases}
\end{align*}
Thus, for $r=0$ we are done.  For $r\geq 1$, suppose that $f_{n,k}(r-1)$  is divisible by $n^{\min\{2,r\}}{2n\choose n}$
for all $n,k$. Applying the relations
$$(n-i)^{2r+1}=(n-i)^{2r-1}(n^2-(2n-i)i)\ \text{and}\ {2n\choose i}(2n-i)i=2n(2n-1){2n-2\choose i-1},$$
we have
\begin{align}
f_{n,k}(r)=n^2 f_{n,k}(r-1)-2n(2n-1)f_{n-1,k-1}(r-1).  \label{eq:rec}
\end{align}
By the induction hypothesis,  $n^2 f_{n,k}(r-1)$ is divisible by
$n^2{2n\choose n}$ and $2n(2n-1)f_{n-1,k-1}(r-1)$ is divisible by
$$2n(2n-1){2n-2\choose n-1}=n^2{2n\choose n}.$$
This proves that $f_{n,k}(r)$ is also divisible by $n^2{2n\choose n}$ for $r\geq 1$.
Finally, it is also easy to check that
$$
n^{-2}{2n\choose n}^{-1}\sum_{k=0}^r f_{n,k}(r)
$$
is an odd integer for $r\geq 1$ by using \eqref{eq:rec} and induction on $r$. The details are left to
the reader.
\qed

Here are some examples for small $r$:
\begin{align}
\sum_{k=1}^n{2n\choose n-k}k^{3} &=\frac{n^2}{2}{2n\choose n},  \nonumber\\
\sum_{k=1}^n{2n\choose n-k}k^{5} &=\frac{n^2}{2}{2n\choose n}(2n-1),\nonumber\\
\sum_{k=1}^n{2n\choose n-k}k^{7} &=\frac{n^2}{2}{2n\choose n}(6n^2-8n+3),\nonumber\\
\sum_{k=1}^n{2n\choose n-k}k^{9} &=\frac{n^2}{2}{2n\choose n}(24n^3-60n^2+54n-17). \nonumber
\end{align}
\begin{rmk}
Note that Shapiro, Woan and Getu \cite{SWG} proved that
\begin{align}
\sum_{n=1}^{\infty}\left(\sum_{k=1}^n k^r B_{n,k}\right)x^n
=\sum_{s=1}^{\lfloor n/2\rfloor}\frac{m(r,s)x^s}{(1-4x)^{(n+1)/2}},  \label{eq:mrs}
\end{align}
where $m(r,s)$ denotes the number of permutations of $\{1,\ldots,r\}$ with $s$ runs and $s$ slides.
Using the binomial theorem, it is easy to give a formula for $\sum_{k=1}^n k^r B_{n,k}$ involving $m(r,s)$
from \eqref{eq:mrs}. A natural question is whether we can deduce Theorem~\ref{thm:half}
from \eqref{eq:mrs}.
\end{rmk}
\section{Proof of Theorem \ref{thm:nnnhalf}}
We will need the Pfaff-Saalsch\"utz identity (see \cite[p.~69]{AAR} or \cite[p.~43, (A)]{CF})
:
\begin{align}
{n_1+n_2\choose n_1+k}{n_2+n_3\choose n_2+k}{n_3+n_1\choose n_3+k}
=\sum_{s=0}^{n_1-k}\frac{(n_1+n_2+n_3-k-s)!}
{s!(s+2k)!(n_1-k-s)!(n_2-k-s)!(n_3-k-s)!},\label{eq:pfaff}
\end{align}
where $\frac{1}{n!}=0$ if $n<0$.

For any positive integers $a_1,\ldots,a_l$, let
$$
C(a_1,\ldots,a_l;k)=\prod_{i=1}^l {a_i+a_{i+1}\choose a_i+k},
$$
where $a_{l+1}=a_1$. Then
\begin{equation}\label{eq:rewriting}
S_{2r+1}(n_1,\ldots,n_{m})
=\frac{(n_1-1)! n_{m}!}{(n_1+n_{m})!}
\sum_{k=1}^{n_1}C(n_1,\ldots,n_{m};k)k^{2r+1}.
\end{equation}

Observe  that  for $m\geq 3$, we have
$$
C(n_1,\ldots,n_m;k)=\frac{(n_2+n_3)!(n_m+n_1)!}{(n_1+n_2)!(n_m+n_3)!}
{n_1+n_2\choose n_1+k}{n_1+n_2\choose n_2+k}C(n_3,\ldots,n_m;k),
$$
and, by letting $n_3\to\infty$ in \eqref{eq:pfaff},
$$
{n_1+n_2\choose n_1+k}{n_1+n_2\choose n_2+k}
=\sum_{s=0}^{n_1-k}\frac{(n_1+n_2)!}{s!(s+2k)!(n_1-k-s)!(n_2-k-s)!}.
$$
Plugging these into \eqref{eq:rewriting} we can write its right-hand side as
\begin{align*}
S_{2r+1}(n_1,\ldots,n_{m})
&=\frac{(n_2+n_3)!(n_1-1)! n_m!}{(n_m+n_3)!}\sum_{k=1}^{n_1}\sum_{s=0}^{n_1-k}
\frac{C(n_3,\ldots,n_{m};k)k^{2r+1}}{s!(s+2k)!(n_1-k-s)!(n_2-k-s)!}  \\
&=\frac{(n_2+n_3)!(n_1-1)! n_m!}{(n_m+n_3)!}\sum_{l=1}^{n_1}\sum_{k=1}^{l}
\frac{C(n_3,\ldots,n_{m};k)k^{2r+1}}{(l-k)!(l+k)!(n_1-l)!(n_2-l)!},
\end{align*}
where $l=s+k$. Now, in the last sum making the substitution
$$
\frac{C(n_3,\ldots, n_{m};k)}{(l-k)!(l+k)!}
=\frac{(n_{m}+n_3)!}{(n_3+l)!(n_{m}+l)!}C(l,n_3,\ldots, n_{m};k),
$$
we obtain the following recurrence relation
\begin{align}
S_{2r+1}(n_1,\ldots,n_{m})
=\sum_{l=1}^{n_1}{n_1-1\choose l-1}{n_2+n_3\choose n_2-l} S_{2r+1}(l,n_3,\ldots,n_{m}). \label{eq:recsr}
\end{align}
By induction on $m$ and using Theorem \ref{thm:half}, we complete the proof.

By iteration of \eqref{eq:recsr} for $r=0,1$, we obtain the following result.
\begin{cor}\label{cor:lambda}
For $m\geq 3$ and all positive integers $n_1,\ldots,n_m$, there hold
\begin{align*}
\sum_{k=1}^{n_1}k\prod_{i=1}^{m} {n_i+n_{i+1}\choose n_i+k}
&=\frac{n_1}{2}{n_1+n_m\choose n_1}
\sum_{\lambda}{\lambda_{m-2}+n_m-1 \choose \lambda_{m-2}}\prod_{i=1}^{m-2}{\lambda_{i-1}-1\choose \lambda_{i}-1}
{n_{i+1}+n_{i+2}\choose n_{i+1}-\lambda_i},  \\
\sum_{k=1}^{n_1}k^3\prod_{i=1}^{m} {n_i+n_{i+1}\choose n_i+k}
&=\frac{n_1 n_m}{2}{n_1+n_m\choose n_1}
\sum_{\lambda}{\lambda_{m-2}+n_m-2 \choose \lambda_{m-2}-1}\prod_{i=1}^{m-2}{\lambda_{i-1}-1\choose \lambda_{i}-1}
{n_{i+1}+n_{i+2}\choose n_{i+1}-\lambda_i},
\end{align*}
where $n_{m+1}=\lambda_0=n_1$ and the sums are over all  sequences
$\lambda=(\lambda_1,\ldots,\lambda_{m-2})$ of positive integers
such that $n_1\geq \lambda_1\geq \cdots \geq \lambda_{m-2}$.
\end{cor}
Note that the following identity was established in \cite{GJZ}:
\begin{equation*}
\sum_{k=-n_1}^{n_1}(-1)^k\prod_{i=1}^m
{n_i+n_{i+1}\choose n_i+k}
={n_1+n_m\choose n_1}\sum_{n_1\geq \lambda_1\geq \cdots \geq \lambda_{m-2}\geq 0}
\prod_{i=1}^{m-2}{\lambda_{i-1}\choose \lambda_i}{n_{i+1}+n_{i+2}
\choose n_{i+1}-\lambda_i},
\end{equation*}
where $n_{m+1}=\lambda_0=n_1$.

\section{Proof of Theorem \ref{thm:ppower}}
We need the following theorem of Lucas (see, for example, \cite{Granville}) and  \cite{SZ} for a recent application.
\begin{thm}[Lucas' theorem]Let $p$ be a prime, and let $a_0,b_0,\ldots, a_m,b_m\in\{0,\ldots,p-1\}$. Then
$$
{\sum_{i=0}^m a_i p^i\choose \sum_{i=0}^m b_i p^i}\equiv \prod_{i=0}^m {a_i\choose b_i}\pmod p.
$$
\end{thm}

Let $r+s\equiv 1\pmod 2$ and $s\geq 1$. Setting $n_1=\cdots=n_m=n$ in Theorem \ref{thm:nnnhalf}, one sees that
$$
\sum_{k=1}^n k^{r+s}{2n\choose n-k}^s
$$
is divisible by ${2n-1\choose n}$. Note that $k{2n\choose n-k}=nB_{n,k}$ is clearly divisible by
$n$. Therefore,
\begin{align*}
\sum_{k=1}^{n}k^{r}B_{n,k}^{s}=n^{-s}\sum_{k=1}^n k^{r+s}{2n\choose n-k}^s
\end{align*}
is divisible by
$$
{2n-1\choose n}\left/\gcd\left({2n-1\choose n},n^s\right)\right..
$$

Now suppose that $n=p^\alpha$ is a prime power.
It follows immediately from Lucas' theorem that
\begin{align*}
{2^{\alpha+1}-1\choose 2^\alpha}={\sum_{k=0}^\alpha 2^k\choose 2^\alpha}\equiv 1\pmod 2,\\
2{2p^\alpha-1\choose p^\alpha}={2p^\alpha\choose p^\alpha}\equiv 2\pmod p\quad(p>2).
\end{align*}
Namely, we always have
$$
{2n-1\choose n}\equiv 1\pmod n,
$$
and thereby
\begin{align}
\gcd\left({2n-1\choose n},n^s\right)=1. \label{eq:2n-1}
\end{align}
This completes the proof.

\medskip
\noindent{\it Remark. } Besides prime powers, there are some other natural numbers satisfying \eqref{eq:2n-1}.
For example, if $p$ and $p^2+p+1$ are both odd primes, then it is easy to see from Lucas' theorem that
$n=p(p^2+p+1)$ satisfies \eqref{eq:2n-1}. In fact, the following are all the natural numbers $n$ less than $500$ not being prime powers
but satisfying \eqref{eq:2n-1}:
$$
39,55,93,111,119,155,161,175,253,275,279,305,317,333,351,363,377,403,407,413,497.
$$
Therefore, Theorem \ref{thm:ppower} is also true for these numbers.

\section{Consequences of Theorem \ref{thm:nnnhalf}}

Letting $n_{2i-1}=m$ and $n_{2i}=n$ for $1\leq i\leq r$ in Theorem~\ref{thm:nnnhalf} and noticing the symmetry of $m$ and $n$,
we obtain
\begin{cor}
For all positive integers $m$, $n$, $r$ and any nonnegative integer $a$,
$$
2\sum_{k=1}^m k^{2a+1} {m+n\choose m+k}^r {m+n\choose n+k}^r
$$
is divisible by $\frac{mn}{\gcd(m,n)}{m+n\choose m}$.
\end{cor}

Letting $n_{3i-2}=l$, $n_{3i-1}=m$ and $n_{3i}=n$ for $1\leq i\leq r$
in Theorem~\ref{thm:nnnhalf}, we obtain
\begin{cor}
For all positive integers $l$, $m$, $n$, $r$ and any nonnegative integer $a$,
$$
2\sum_{k=1}^l k^{2a+1}{l+m\choose l+k}^r{m+n\choose m+k}^r {n+l\choose n+k}^r
$$
is divisible by $\frac{lm}{\gcd(l,m)}{l+m\choose l}$, $\frac{mn}{\gcd(m,n)}{m+n\choose m}$ and $\frac{nl}{\gcd(n,l)}{n+l\choose n}$.
\end{cor}

Letting $m=2r+s$, $n_1=n_3=\cdots=n_{2r-1}=n+1$ and letting all the other $n_i$ be $n$ in
Theorem~\ref{thm:nnnhalf}, we get
\begin{cor}\label{cor:nn+1}
For all positive integers $r$, $s$, $n$ and any nonnegative integer $a$,
\begin{align*}
\sum_{k=1}^n k^{2a+1} {2n+1\choose n+k+1}^r{2n+1\choose n+k}^r{2n\choose n+k}^s
\end{align*}
is divisible by $\frac{n(n+1)}{2}{2n+1\choose n}$.
\end{cor}

Clearly Theorem~\ref{thm:nnnhalf} can be restated in the following form.

\begin{thm}\label{thm:rennnhalf}
For all positive integers $n_1,\ldots,n_m$ and any nonnegative integer $r$,
$$
2(n_1-1)!\prod_{i=1}^m\frac{(n_i+n_{i+1})!}{(2n_i)!}
\sum_{k=1}^{n_1} k^{2r+1}\prod_{i=1}^m {2n_i\choose n_i+k},
$$
where $n_{m+1}=0$, is an integer.
\end{thm}

It is not hard to see that, for all positive integers $m$ and $n$, the expression
$\frac{(2m)!(2n)!}{2(m+n)!m!n!}$ is an integer by considering
the power of a prime dividing a factorial. Letting $n_1=\cdots= n_r=m$ and
$n_{r+1}=\cdots=n_{r+s}=n$ in Theorem~\ref{thm:rennnhalf} and noticing the symmetry of $m$ and $n$, we obtain
\begin{cor}\label{cor:mn}
For all positive integers $m$, $n$, $r$, $s$ and any nonnegative integer $a$,
$$
\sum_{k=1}^m k^{2a+1}{2m\choose m+k}^r{2n\choose n+k}^s
$$
is divisible by $\frac{(2m)!(2n)!}{2(m+n)!(m-1)!(n-1)!\gcd(m,n)}$.
\end{cor}

In particular, we find that
$$
\sum_{k=1}^n k^{2a+1}{4n\choose 2n+k}^r{2n\choose n+k}^s
$$
is divisible by $n{4n\choose n}$, and
$$
\sum_{k=1}^n k^{2a+1}{6n\choose 3n+k}^r{2n\choose n+k}^s
$$
is divisible by
$\frac{(6n)!(2n-1)!}{(4n)!(3n-1)!(n-1)!}$.

Using Lucas' theorem, similarly to Theorem \ref{thm:ppower}, we can deduce the following result immediately.
\begin{cor}\label{cor:n2n}
Let $n$ be a power of $2$. Let $r\geq 0$ and $s,t\geq 1$ such that $r+s+t\equiv 1\pmod 2$. Then
\begin{align*}
\sum_{k=1}^{n}k^{r}B_{2n,k}^{s}B_{n,k}^t
\end{align*}
is divisible by ${4n-1\choose n-1}$.
\end{cor}

{}From Theorem~\ref{thm:rennnhalf} it is easy to see that
$$
2(n_1-1)!\prod_{i=1}^m\frac{(n_i+n_{i+1})!}{(2n_i)!}
\sum_{k=1}^{n_1} k^{2a+1}\prod_{i=1}^m {2n_i\choose n_i+k}^{r_i},
$$
where $n_{m+1}=0$, is a nonnegative integer for all $r_1,\ldots,r_m\geq 1$.
For $m=3$, letting $(n_1,n_2,n_3)$ be $(n,3n,2n)$, $(2n,n,3n)$, or
$(2n,n,4n)$, and noticing the symmetry of $n_1$ and $n_3$, we obtain the following two corollaries.
\begin{cor}\label{cor:rst-246n}
For all positive integers $r$, $s$, $t$, $n$ and any nonnegative integer $a$,
$$
\sum_{k=1}^n k^{2a+1} {6n\choose 3n+k}^r {4n\choose 2n+k}^s
{2n\choose n+k}^t
$$
is divisible by both $n{6n\choose n}$ and $3n{6n\choose 3n}$.
\end{cor}

\begin{cor}\label{cor:rst-248n}
For all positive integers $r$, $s$, $t$, $n$ and any nonnegative integer $a$,
$$
\sum_{k=1}^n k^{2a+1}{8n\choose 4n+k}^r {4n\choose 2n+k}^s
{2n\choose n+k}^t
$$
is divisible by $2n{8n\choose 3n}$.
\end{cor}

\section{Concluding remarks and open problems}

It is easy to prove that
$$
\sum_{k=1}^n{2n\choose n-k}k^{2}=4^{n-1} n
$$
(see \cite{SWG}). Furthermore, similarly to \eqref{eq:sumnkl}, there holds
\begin{align}
\sum_{k=1}^n{2n\choose n-k}k^{2r}
=\sum_{k=0}^{r-1} 4^{n-k-1}
\sum_{i=0}^k{2n\choose i}{2k-2n\choose k-i}(n-i)^{2r-1},\quad r\geq 1.  \label{eq:k2r}
\end{align}
In particular, for $r=2,3,4,5$, we have
\begin{align*}
\sum_{k=1}^n{2n\choose n-k}k^{4} &=2^{2n-3}n(3n-1),    \\
\sum_{k=1}^n{2n\choose n-k}k^{6} &=2^{2n-4}n(15n^2-15n+4),    \\
\sum_{k=1}^n{2n\choose n-k}k^{8} &=2^{2n-5}n(105n^3-210n^2+147n-34),\\
\sum_{k=1}^n{2n\choose n-k}k^{10}&=2^{2n-6}n(945n^4-3150n^3+4095n^2-2370n+496).
\end{align*}

Let $\alpha(N)$ denote the number of $1$'s in the binary expansion of $N$. For example, $\alpha(255)=8$ and $\alpha(256)=1$.
\begin{conj}\label{conj:alpha}
Let $n,r\geq 1$. Then
\begin{align}
\sum_{k=1}^n{2n\choose n-k}k^{2r}  \label{eq:sumk2r}
\end{align}
is divisible by $2^{2n-\min\{\alpha(n),\,\alpha(r)\}-1}$.
\end{conj}

Note that it is clear from \eqref{eq:k2r} that \eqref{eq:sumk2r} is divisible by $4^{n-r}$, assuming that $n\geq r$.
The statistic $\alpha(n)$ appears in Conjecture \ref{conj:alpha} because the divisibility of \eqref{eq:sumk2r}
also depends on $n$.

\begin{conj}Let $n,r\geq 1$. Then
$$
\sum_{k=1}^n B_{n,k}^{2r+1}\equiv {2n-1\choose n}\mod {{2n\choose n}}
$$
if and only if $n= 2^a-2^b$ for some $0\leq b<a$.
\end{conj}

For example, we have
\begin{align*}
&\sum_{k=1}^{7} B_{7,k}^{3}=354331692\equiv 1716\pmod{3432}, \\
&\sum_{k=1}^{12} B_{12,k}^{3}=96906387191038334\equiv 1352078\pmod{2704156}, \\
&\sum_{k=1}^{13} B_{13,k}^{3}=5066711735118128200\equiv 0\pmod{10400600},\\
&\sum_{k=1}^{16} B_{16,k}^{3}=786729115199980286001225 \equiv 0\pmod{601080390}.
\end{align*}
\begin{conj}
Let $r\geq 0$ and $s\geq 1$. Then for $n\geq 4^s-1$,
\begin{align*}
\sum_{k=1}^{n}k^{2r+1}B_{n,k}^{2s}
\equiv
\begin{cases}
\displaystyle {2n-1\choose n}4^{s-1},&\text{if $n=4^s-1$ or $n=2^a+1$,}\\[5pt]
0,&\text{otherwise.}
\end{cases}
\mod{{2n\choose n}4^{s-1}}.
\end{align*}
\end{conj}

\begin{conj}
Let $n,r,s\geq 1$. Then
\begin{align*}
\sum_{k=1}^{n}k^{2r}B_{n,k}^{2s+1}
\equiv
\begin{cases}
\displaystyle {2n-1\choose n},&\text{if $n=2^a-1$,}\\[5pt]
0,&\text{otherwise.}
\end{cases}
\mod{{2n\choose n}}.
\end{align*}
\end{conj}

The following statement is a generalization of Corollary \ref{cor:n2n}.
\begin{conj}Let $m, n, s, t\geq 1$, and $r\geq 0$ such that $r+s+t\equiv 1\pmod 2$. Then
$$
\sum_{k=1}^n k^r B_{m,k}^{s}B_{n,k}^{t}
$$
is divisible by $\frac{1}{2}\frac{(2m)!(2n)!}{m!n!(m+n)!}$.
\end{conj}

Furthermore, for some special cases, we have
\begin{conj}Let $n,r,s\geq 1$ such that $r\not\equiv s\pmod 2$. Then
$$
\sum_{k=1}^n B_{n,k}^{r}B_{2n,k}^{s}
\equiv
\begin{cases}
\displaystyle {4n\choose n}2^{\min\{2r,s\}-2},&\text{if $n=\frac{2^a(2^{2b+1}+1)}{3}$,}\\[5pt]
0,&\text{otherwise.}
\end{cases}
\mod{{4n\choose n}2^{\min\{2r,s\}-1}}.
$$
\end{conj}

\begin{conj}Let $n,r,s\geq 1$ such that $r\not\equiv s\pmod 2$. Then
\begin{align*}
\sum_{k=1}^n B_{n,k}^{r}B_{n+1,k}^{s}
\equiv
\begin{cases}
\displaystyle {2n\choose n}2^{\min\{r,s\}-1},&\!\!\!\text{if\ $\begin{cases} r<s, \ n=2^a-2^b\ {\rm(}b\geq 1{\rm)}
\\ r>s,\ n=2^a-1 \end{cases}$}\\[5pt]
0,&\text{otherwise.}
\end{cases}
\!\!\!\!\bmod{{2n\choose n}2^{\min\{r,s\}}},
\end{align*}
\end{conj}

The following two conjectures are refinements of Corollaries \ref{cor:rst-246n} and \ref{cor:rst-248n}
for $2a+1=r+s+t$.

\begin{conj}For any positive integers $r,s,t$ such that $r+s+t\equiv1\pmod 2$,
$$
\sum_{k=1}^n B_{n,k}^{r}B_{2n,k}^{s}B_{3n,k}^{t}
$$
is divisible by both $\frac{1}{3}{6n\choose n}$ and ${6n\choose 3n}$.
\end{conj}

\begin{conj}For any positive integers $r,s,t$ such that $r+s+t\equiv1\pmod 2$,
$$
\sum_{k=1}^n B_{n,k}^{r}B_{2n,k}^{s}B_{4n,k}^{t}
$$
is divisible by ${8n\choose 3n}$.
\end{conj}

Finally it is natural to introduce the $q$-Catalan triangles with entries given by
$$
B_{n,k}(q):=\frac{1-q^k}{1-q^n}{2n\brack n-k},\quad 1\leq k\leq n,
$$
where the $q$-binomial coefficient is defined by
$$
{M\brack N}:=
\begin{cases}
\displaystyle\prod_{i=1}^N\frac{1-q^{M-N+i}}{1-q^i},&\text{if $0\leq N\leq M$,}\\[5pt]
0,&\text{otherwise.}
\end{cases}
$$
Let $\Phi_n(x)$ be the $n$-th cyclotomic polynomial. It is not difficult to see from \cite[Eq.~(10)]{KW}
or \cite[Proposition 2.2]{GZ06} that
$$
B_{n,k}(q)=\left(\prod_i\Phi_i(q)\right)\left(\prod_{d}\Phi_d(q)\right),
$$
where the first product is over all positive integers $i$ such that $i\mid k$ and $i\nmid n$ and
the second product is over all positive integers $d$ such that
$\lfloor (n-k)/d\rfloor+\lfloor (n+k)/d\rfloor<\lfloor 2n/d\rfloor$  and $d$ does not divide $n$.


Now, the obvious identity
$$
\frac{1-q^k}{1-q^n}{2n\brack n-k}q^{k\choose 2}
={2n-1\brack n-k}q^{k\choose 2}-{2n-1\brack n-k-1}q^{k+1\choose 2}
$$
implies the following  $q$-analogue of \eqref{eq:shap}:
\begin{align}
\sum_{k=1}^n\frac{1-q^k}{1-q^n}{2n\brack n-k}q^{k\choose 2}=\frac{1}{1+q^n}{2n\brack n}.  \label{eq:qshap}
\end{align}
In view  of the identity \eqref{eq:qshap} and the results in \cite{GJZ},  it would be interesting
to find a $q$-analogue of Theorems {\rm\ref{thm:half}} and {\rm\ref{thm:nnnhalf}}.

\medskip
\noindent{\bf Acknowledgements.} This work was supported by the projects MIRA 2007 and MIRA 2008 of R\'egion Rh\^one-Alpes and
was mainly done during the first author's visit to
Institut Camille Jordan of Univerit\'e
Lyon I. The first author was also sponsored by Shanghai Educational Development Foundation under the Chenguang
Project (\#2007CG29), Shanghai Leading Academic Discipline Project
(\#B407), and the National Science Foundation of China (\#10801054).

\renewcommand{\baselinestretch}{1}

\end{document}